\documentclass[11pt]{tran-l}
\usepackage[ansinew]{inputenc}
\usepackage[centertags]{amsmath}
\usepackage{exscale}
\usepackage{relsize}
\usepackage{amsfonts}
\usepackage{amssymb}
\usepackage{amsthm}
\usepackage{newlfont}

\mathchardef\ordinarycolon\mathcode`\:
  \mathcode`\:=\string"8000
  \begingroup \catcode`\:=\active
    \gdef:{\mathrel{\mathop\ordinarycolon}}
  \endgroup
\def\N{{\Bbb N}}

\def\R{{\Bbb R}}
\newtheorem*{mthm*}{Main Theorem}
\newtheorem{mthm}{Main Theorem}
\newtheorem{thm}{Theorem}
\newtheorem{lem}[thm]{Lemma}
\newtheorem{prop}[thm]{Proposition}

\newtheorem{rem}[thm]{Remark}
\newtheorem{Cor}[thm]{Corollary}

\newcommand{\D}{\mathcal{D}}
\newcommand{\B}{\mathcal{B}}

\newcommand{\aut}{\emph{Aut}}
\newcommand{\Aut}{\mbox{Aut}}

\newcommand{\G}{\mathit{\Gamma}}
\newcommand{\Si}{\mathit{\Sigma}}
\newcommand{\q}{\overline{q}}

\begin{document}
\title[Highly Symmetric Combinatorial Designs]
{A Census of Highly Symmetric\\ Combinatorial Designs}

\author{Michael Huber}

\address{Mathematisches Institut der Universit\"{a}t T\"{u}bingen, Auf der Morgenstelle~10,
D-72076~T\"{u}bingen, Germany}

\email{michael.huber@uni-tuebingen.de}

\subjclass[2000]{Primary 51E10; Secondary 05B05, 20B25}

\keywords{Steiner designs, flag-transitive group of automorphisms,
\mbox{$3$-homogeneous} permutation groups}


\date{January 3, 2006; and in revised form January 29, 2007}

\commby{}


\begin{abstract}
As a consequence of the classification of the finite simple groups,
it has been possible in recent years to characterize Steiner
\mbox{$t$-designs}, that is $t$-$(v,k,1)$ designs, mainly for
\mbox{$t=2$}, admitting groups of automorphisms with sufficiently
strong symmetry properties. However, despite the finite simple group
classification, for Steiner \mbox{$t$-designs} with $t>2$ most of
these characterizations have remained long-standing challenging
problems. Especially, the determination of all flag-transitive
Steiner \mbox{$t$-designs} with $3 \leq t\leq 6$ is of particular
interest and has been open for about 40 years
(cf.~\cite[p.\,147]{Del1992} and~\cite[p.\,273]{Del1995}, but
presumably dating back to 1965).

\noindent The present paper continues the author's
work~\cite{Hu2001,Hu2005,Hu2006} of classifying all flag-transitive
Steiner \mbox{$3$-designs} and \mbox{$4$-designs}. We give a
complete classification of all flag-transitive Steiner
\mbox{$5$-designs} and prove furthermore that there are no
non-trivial flag-transitive Steiner \mbox{$6$-designs}. Both results
rely on the classification of the finite \mbox{$3$-homogeneous}
permutation groups. Moreover, we survey some of the most general
results on highly symmetric Steiner \mbox{$t$-designs}.
\end{abstract}

\maketitle

\section{Introduction}\label{intro}

For positive integers $t \leq k \leq v$ and $\lambda$, we define a
\mbox{\emph{$t$-$(v,k,\lambda)$ design}} to be a finite incidence
structure \mbox{$\D=(X,\B,I)$}, where $X$ denotes a set of
\emph{points}, $\left| X \right| =v$, and $\B$ a set of
\emph{blocks}, $\left| \B \right| =b$, with the properties that each
block $B \in \B$ is incident with $k$ points, and each
\mbox{$t$-subset} of $X$ is incident with $\lambda$ blocks. A
\emph{flag} of $\D$ is an incident point-block pair $(x,B) \in I$
with $x \in X$ and $B \in \B$. We consider automorphisms of $\D$ as
pairs of permutations on $X$ and $\B$ which preserve incidence, and
call a group \mbox{$G \leq \Aut (\D)$} of automorphisms of $\D$
\emph{flag-transitive} (respectively \emph{block-transitive},
\emph{point \mbox{$t$-transitive}}, \emph{point
\mbox{$t$-homogeneous}}) if $G$ acts transitively on the flags
(respectively transitively on the blocks, $t$-transitively on the
points, $t$-homogeneously on the points) of $\D$. For short, $\D$ is
said to be, e.g., flag-transitive if $\D$ admits a flag-transitive
group of automorphisms. For historical reasons, a
\mbox{$t$-$(v,k,\lambda)$ design} with \mbox{$\lambda =1$} is called
a \emph{Steiner \mbox{$t$-design}} (sometimes also known as a
\emph{Steiner system}). We note that in this case each block is
determined by the set of points which are incident with it, and thus
can be identified with a $k$-subset of $X$ in a unique way. If
$t<k<v$ holds, then we speak of a \emph{non-trivial} Steiner
\mbox{$t$-design}.

As a consequence of the classification of the finite simple groups,
it has been possible in recent years to characterize Steiner
\mbox{$t$-designs}, mainly for \mbox{$t=2$}, admitting groups of
automorphisms with sufficiently strong symmetry properties. However,
despite the classification of the finite simple groups, for Steiner
\mbox{$t$-designs} with $t>2$ most of these characterizations have
remained long-standing challenging problems. Especially, the
determination of all flag-transitive Steiner \mbox{$t$-designs} with
$3 \leq t\leq 6$ is of particular interest and has been open for
about 40 years (cf.~\cite[p.\,147]{Del1992}
and~\cite[p.\,273]{Del1995}, but presumably dating back to 1965).

The present paper continues the author's
work~\cite{Hu2001,Hu2005,Hu2006} of classifying all flag-transitive
Steiner \mbox{$3$-designs} and \mbox{$4$-designs}. We give a
complete classification of all flag-transitive Steiner
\mbox{$5$-designs} in Section~\ref{flag5designs} and prove
furthermore in Section~\ref{flag6designs} that there are no
non-trivial flag-transitive Steiner \mbox{$6$-designs}. Both results
rely on the classification of the finite \mbox{$3$-homogeneous}
permutation groups, which itself depends on the finite simple group
classification. Summarizing our results in this paper, we state:

\medskip

The classification of all non-trivial Steiner $t$-designs with $t=5$
or $6$ admitting a flag-transitive group of automorphisms is as
follows.

\begin{mthm*}\label{mainthm}
Let $\D=(X,\B,I)$ be a non-trivial Steiner \mbox{$t$-design} with
$t=5$ or $6$. Then \mbox{$G \leq \aut(\D)$} acts flag-transitively
on $\D$ if and only if one of the following occurs:

\begin{enumerate}
\item[(1)] $\D$ is isomorphic to the Witt \mbox{$5$-$(12,6,1)$}
design, and \mbox{$G \cong M_{12}$},

\medskip

\item[(2)] $\D$ is isomorphic to the Witt \mbox{$5$-$(24,8,1)$}
design, and \mbox{$G \cong PSL(2,23)$} or \mbox{$G \cong M_{24}$}.
\end{enumerate}
\end{mthm*}

Referring to the author's work mentioned above, we present the
complete determination of all flag-transitive Steiner
\mbox{$t$-designs} with $t \geq 3$ in Section~\ref{census}.
Moreover, we give in this context a survey on some of the most
general results on highly symmetric Steiner \mbox{$t$-designs}.

\bigskip


\section{Classifications of Highly Symmetric Combinatorial Designs}\label{census}

In the sequel, we survey classification results of  highly symmetric
Steiner \mbox{$t$-designs}. For detailed descriptions of the
respective designs and their groups of automorphisms as well as for
further surveys concerning in particular highly symmetric Steiner
\mbox{$2$-designs}, we refer
to~\cite[Sect.\,1,\,2]{Buek1988},~\cite[Ch.\,2.3,\,2.4,\,4.4]{Demb1968},
~\cite{Kant1985b} and~\cite{Kant1993}.

As presumably one of the first most general results, all point
\mbox{$2$-transitive} Steiner \mbox{$2$-designs} were characterized
by W.~M.~Kantor~\cite[Thm.\,1]{Kant1985}, using the classification
of the finite \mbox{$2$-transitive} permutation groups.

\begin{thm}{\em (Kantor~1985).}\label{2trsdes-A}
Let $\D=(X,\B,I)$ be a non-trivial Steiner \mbox{$2$-design}, and let
\mbox{$G \leq \aut(\D)$} act point \mbox{$2$-transitively} on $\D$.
Then one of the following holds:

\smallskip

\begin{enumerate}

\item[(1)] ${\D}$ is isomorphic to the
$2$-$(\frac{q^d-1}{q-1},q+1,1)$ design whose points and blocks are
the points and lines of the projective space $PG(d-1,q)$, and
$PSL(d,q) \leq G \leq P \mathit{\Gamma} L(d,q)$, or $(d-1,q)=(3,2)$
and $G \cong A_7$,

\medskip

\item[(2)] $\D$ is isomorphic to a Hermitian unital $U_H(q)$ of order $q$,
and  \linebreak $PSU(3,q^2) \leq G \leq P \mathit{\Gamma} U
(3,q^2)$,

\medskip

\item[(3)] $\D$ is isomorphic to a Ree unital $U_R(q)$ of order $q$ with
$q=3^{2e+1}>3$, and  $Re(q) \leq G \leq \aut(Re(q))$,

\medskip

\item[(4)] $\D$ is isomorphic to the $2$-$(q^d,q,1)$ design
whose points and blocks are the points and lines of the affine space
$AG(d,q)$, and one of the following holds (where $G_0$ denotes the
stabilizer of $0 \in X$):

\smallskip

\begin{enumerate}

\item[(i)] $G \leq A \mathit{\Gamma} L(1,q^d)$,

\smallskip

\item[(ii)] $G_0 \unrhd SL(\frac{d}{a},q^a)$, $d \geq 2a$,

\smallskip

\item[(iii)] $G_0 \unrhd Sp(\frac{2d}{a},q^a)$, $d \geq 2a$,

\smallskip

\item[(iv)] $G_0 \unrhd G_2(q^a)'$, $q$ even, $d=6a$,

\smallskip

\item[(v)] $G_0 \unrhd SL(2,3)$ or $SL(2,5)$, $v=q^2$,
$q=5,7,9,11,19,23,29$ or $59$,

\smallskip

\item[(vi)] $G_0 \unrhd SL(2,5)$, or $G_0$ contains a normal extraspecial subgroup $E$ of
order $2^5$ and $G_0/E$ is isomorphic to a subgroup of $S_5$,
$v=3^4$,

\smallskip

\item[(vii)] $G_0 \cong SL(2,13)$, $v=3^6$,

\end{enumerate}

\medskip

\item[(5)] $\D$ is isomorphic to the affine nearfield plane
$A_9$ of order $9$, and $G_0$ as in $\emph{(4)(vi)}$,

\medskip

\item[(6)] $\D$ is isomorphic to the affine Hering plane
$A_{27}$ of order $27$, and $G_0$ as in $\emph{(4)(vii)}$,

\medskip

\item[(7)] $\D$ is isomorphic to one of the two Hering spaces
$2$-$(9^3,9,1)$, and $G_0$ as in $\emph{(4)(vii)}$.

\end{enumerate}
\end{thm}

As an easy implication, W.~M.~Kantor~\cite[Thm.\,3]{Kant1985}
obtained moreover the classification of all point $t$-transitive
Steiner \mbox{$t$-designs} with \mbox{$t > 2$}.

\smallskip

Certainly, among the highly symmetric properties of incidence
structures, flag-transitivity is a particularly important and
natural one. Even long before the aforementioned classification of
the finite simple groups, a general study of flag-transitive Steiner
\mbox{$2$-designs} was introduced by D.~G.~Higman and
J.~E.~McLaughlin~\cite{HigMcL1961} proving that a flag-transitive
group $G \leq \Aut(\D)$ of automorphisms of a Steiner
\mbox{$2$-design} $\D$ is necessarily primitive on the points of
$\D$. They posed the problem of classifying all finite
flag-transitive projective planes, and showed that such planes are
Desarguesian if its orders are suitably restricted. Much later
W.~M.~Kantor~\cite{Kant1987} determined all such planes apart from
the still open case when the group of automorphisms is a Frobenius
group of prime degree. His proof involves detailed knowledge of
primitive permutation groups of odd degree based on the
classification of the finite simple groups. In a big common effort,
F.~Buekenhout, A.~Delandtsheer, J.~Doyen, P.~B.~Kleidman,
M.~W.~Liebeck, and
J.~Saxl~\cite{Buek1990,Del2001,Kleid1990,Lieb1998,Saxl2002}
essentially characterized all finite flag-transitive linear spaces,
that is flag-transitive Steiner \mbox{$2$-designs}. Their result,
which also relies on the finite simple group classification, starts
with the result of Higman and McLaughlin and uses the O'Nan-Scott
Theorem for finite primitive permutation groups. For the incomplete
case with a $1$-dimensional affine group of automorphisms, we refer
to~\cite[Sect.\,4]{Buek1990} and~\cite[Sect.\,3]{Kant1993}.

\begin{thm}{\em (Buekenhout et al.~1990).}\label{flag2des}
Let ${\D}=(X,\B,I)$ be a Steiner \mbox{$2$-design}, and let \mbox{$G \leq
\aut(\D)$} act flag-transitively on $\D$. Then one of the following
occurs:

\smallskip

\begin{enumerate}

\item[(1)] ${\D}$ is isomorphic to the $2$-$(q^d,q,1)$ design
whose points and blocks are the points and lines of the affine space
$AG(d,q)$, and one of the following holds:

\smallskip

\begin{enumerate}

\item[(i)] $G$ is $2$-transitive (hence as in~\emph{Theorem~\ref{2trsdes-A} (4)}),

\smallskip

\item[(ii)] $d=2$, $q=11$ or $23$, and $G$ is one of the three
solvable flag-transitive groups given
in~\cite[Table\,II]{Foul1964a},

\smallskip

\item[(iii)] $d=2$, $q=9,11,19,29$ or $59$,
$G_0^{(\infty)} \cong SL(2,5)$ (where $G_0^{(\infty)}$ denotes the
last term in the derived series of $G_0$), and $G$ is given
in~\cite[Table\,II]{Foul1964a},

\smallskip

\item[(iv)] $d=4$, $q=3$, and $G_0 \cong SL(2,5)$,

\end{enumerate}

\medskip

\item[(2)] $\D$ is isomorphic to a non-Desarguesian affine
translation plane. More precisely, one of the following holds:

\smallskip

\begin{enumerate}

\item[(i)] $\D$ is isomorphic to a L\"{u}neburg-Tits plane
\emph{Lue}$(q^2)$ of order $q^2$ with $q=2^{2e+1}>2$, and $Sz(q)
\leq G_0 \leq \aut(Sz(q))$,

\smallskip

\item[(ii)] $\D$ is isomorphic to the affine Hering plane
$A_{27}$ of order $27$, and $G_0 \cong SL(2,13)$,

\smallskip

\item[(iii)] $\D$ is isomorphic to the affine nearfield plane
$A_9$ of order $9$, and $G$ is one of the seven flag-transitive
subgroups of $\aut(A_9)$, described in~\cite[§\,5]{Foul1964b},

\end{enumerate}

\medskip

\item[(3)] $\D$ is isomorphic to one of the two Hering spaces
$2$-$(9^3,9,1)$, \linebreak and $G_0 \cong  SL(2,13)$,

\medskip

\item[(4)] $\D$ is isomorphic to the
$2$-$(\frac{q^d-1}{q-1},q+1,1)$ design whose points and blocks are
the points and lines of the projective space $PG(d-1,q)$, and
$PSL(d,q) \leq G \leq P \mathit{\Gamma} L(d,q)$, or $(d-1,q)=(3,2)$
and $G \cong A_7$,

\medskip

\item[(5)] $\D$ is isomorphic to a Hermitian unital $U_H(q)$ of order $q$,
\linebreak and $PSU(3,q^2) \leq G \leq P \mathit{\Gamma} U (3,q^2)$,

\medskip

\item[(6)] $\D$ is isomorphic to a Ree unital $U_R(q)$ of order $q$ with
$q=3^{2e+1}>3$, and $Re(q) \leq G \leq \aut(Re(q))$,

\medskip

\item[(7)] $\D$ is isomorphic to a Witt-Bose-Shrikhande space
$W(q)$ with \linebreak \mbox{$q=2^d \geq 8$}, and $PSL(2,q) \leq G
\leq P \mathit{\Gamma} L(2,q)$,

\medskip

\item[(8)] $G \leq A \mathit{\Gamma} L (1,q)$.

\end{enumerate}
\end{thm}

\medskip

Investigating \mbox{$t$-designs} $\D$ for arbitrary $\lambda$, but
large $t$, P.~J.~Cameron and
C.~E.~Praeger~\cite[Thm.\,1.1\,and\,2.1]{CamPrae1993} showed that
for $t \geq 7$ the flag-transitivity, respectively for $t \geq 8$
the block-transitivity of \mbox{$G \leq \Aut (\D)$} implies at least
its point \mbox{$4$-homogeneity} and proved the following result:

\begin{thm}{\em (Cameron and Praeger~1993).}\label{flag7des}
Let $\D=(X,\B,I)$ be a \mbox{$t$-$(v,k,\lambda)$} design. If
\mbox{$G \leq \aut(\D)$} acts block-transitively on $\D$, then $t
\leq 7$, while if \mbox{$G \leq \aut(\D)$} acts flag-transitively on
$\D$, then $t \leq 6$.
\end{thm}

However, especially the determination of all flag-transitive Steiner
\linebreak \mbox{$t$-designs} with $3 \leq t \leq 6$ has remained of
particular interest, and even the classification of all
flag-transitive Steiner \mbox{$3$-designs} has been known as ''a
long-standing and still open problem'' (cf.~\cite[p.\,147]{Del1992}
and~\cite[p.\,273]{Del1995}). Presumably,
H.~L\"{u}neburg~\cite{Luene1965} in 1965 has been the first dealing
with part of this problem characterizing flag-transitive Steiner
quadruple systems (i.e., Steiner \mbox{$3$-designs} with block size
$k=4$) under the additional strong assumption that every
non-identity element of the group of automorphisms fixes at most two
distinct points. This result has been generalized in 2001 by the
author~\cite{Hu2001}, omitting the additional assumption on the
number of fixed points. Recently, the author~\cite{Hu2005,Hu2006}
completely determined all flag-transitive Steiner \mbox{$3$-designs}
and \mbox{$4$-designs} using the classification of the finite
\mbox{$2$-transitive} permutation groups. In the present paper, the
remaining investigations of all flag-transitive Steiner
\mbox{$5$-designs} and \mbox{$6$-designs} are given, utilizing the
classification of the finite \mbox{$3$-homogeneous} permutation
groups. Summarizing the author's results, the complete determination
of all non-trivial Steiner \mbox{$t$-designs} with $t \geq 3$
admitting a flag-transitive group of automorphisms can now be stated
as follows.

\begin{thm}{\em (Huber~2005/06).}\label{flagtdes}
Let $\D=(X,\B,I)$ be a non-trivial Steiner \mbox{$t$-design} with
\mbox{$t \geq 3$}. Then \mbox{$G \leq \aut(\D)$} acts
flag-transitively on $\D$ if and only if one of the following
occurs:

\smallskip

\begin{enumerate}

\item[(1)] $\D$ is isomorphic to the \mbox{$3$-$(2^d,4,1)$} design
whose points and blocks are the points and planes of the affine
space $AG(d,2)$, and one of the following holds:

\smallskip

\begin{enumerate}
\item[(i)] $d \geq 3$, and $G \cong AGL(d,2)$,

\smallskip

\item[(ii)] $d=3$, and $G \cong AGL(1,8)$ or $A \mathit{\Gamma}
L(1,8)$,

\smallskip

\item[(iii)] $d=4$, and $G_0 \cong A_7$,

\smallskip

\item[(iv)] $d=5$, and $G \cong A \mathit{\Gamma} L(1,32)$,
\end{enumerate}

\medskip

\item[(2)] $\D$ is isomorphic to a \mbox{$3$-$(q^e +1,q+1,1)$}
design whose points are the elements of the projective line
\mbox{$GF(q^e) \cup \{\infty\}$} and whose blocks \linebreak are the
images of \mbox{$GF(q) \cup \{\infty\}$} under $PGL(2,q^e)$
(respectively \linebreak $PSL(2,q^e)$, $e$ odd) with a prime power
\mbox{$q \geq 3$}, \mbox{$e \geq 2$}, and the derived design at any
given point is isomorphic to the \mbox{$2$-$(q^e,q,1)$} design whose
points and blocks are the points and lines of $AG(e,q)$, and
\mbox{$PSL(2,q^e) \leq G \leq P \mathit{\Gamma} L (2,q^e)$},

\medskip

\item[(3)] $\D$ is isomorphic to a \mbox{$3$-$(q+1,4,1)$} design
whose points are the elements of \mbox{$GF(q) \cup \{\infty\}$} with
a prime power \mbox{$q \equiv 7$ $($\emph{mod} $12)$} and whose
blocks are the images of \mbox{$\{0,1,\varepsilon,\infty\}$} under
$PSL(2,q)$, where $\varepsilon$ is a primitive sixth root of unity
in $GF(q)$, and the derived design at any given point is isomorphic
to the Netto triple system $N(q)$, and \mbox{$PSL(2,q) \leq G \leq P
\mathit{\Sigma} L (2,q)$},

\medskip

\item[(4)] $\D$ is isomorphic to one of the following Witt designs:

\smallskip

\begin{enumerate}

\item[(i)] the \mbox{$3$-$(22,6,1)$} design, and \mbox{$G \unrhd M_{22}$},

\smallskip

\item[(ii)] the \mbox{$4$-$(11,5,1)$} design, and \mbox{$G \cong M_{11}$},

\smallskip

\item[(iii)] the \mbox{$4$-$(23,7,1)$} design, and \mbox{$G \cong M_{23}$},

\smallskip

\item[(iv)] the \mbox{$5$-$(12,6,1)$} design, and \mbox{$G \cong M_{12}$},

\smallskip

\item[(v)] the \mbox{$5$-$(24,8,1)$} design, and \mbox{$G \cong PSL(2,23)$} or \mbox{$G \cong M_{24}$}.
\end{enumerate}
\end{enumerate}
\end{thm}

\medskip

We remark that the Steiner \mbox{$3$-designs} in Part (1) (ii) with
\mbox{$G \cong AGL(1,8)$} and (iv) with \mbox{$G \cong A
\mathit{\Gamma} L(1,32)$} as well as the Steiner \mbox{$5$-design}
in Part (4) with \mbox{$G \cong PSL(2,23)$} are sharply
flag-transitive, and furthermore, concerning Part (4)~(v), that
$M_{24}$ as the full group of automorphisms of $\D$ contains only
one conjugacy class of subgroups isomorphic to $PSL(2,23)$.

\bigskip


\section{Definitions and Preliminary Results}\label{Prelim}

If $\D=(X,\B,I)$ is a \mbox{$t$-$(v,k,\lambda)$} design with $t \geq
2$, and $x \in X$ arbitrary, then the \emph{derived} design with
respect to $x$ is \mbox{$\D_x=(X_x,\B_x, I_x)$}, where $X_x = X
\backslash \{x\}$, \mbox{$\B_x=\{B \in \B: (x,B)\in I\}$} and $I_x=
I \!\!\mid _{X_x \times \; \B_x}$. In this case, $\D$ is also called
an \emph{extension} of $\D_x$. Obviously, $\D_x$ is a
\mbox{$(t-1)$-$(v-1,k-1,\lambda)$} design.

Let $G$ be a permutation group on a non-empty set $X$. We call $G$
\emph{semi-regular} if the identity is the only element that fixes
any point of $X$. If additionally $G$ is transitive, then it is said
to be \emph{regular}. Furthermore, for $x \in X$, the orbit $x^G$
containing $x$ is called \emph{regular} if it has length $\left| G
\right|$. If \mbox{$\{x_1,\ldots,x_m\} \subseteq X$}, let
$G_{\{x_1,\ldots,x_m\}}$ be its setwise stabilizer and
$G_{x_1,\ldots,x_m}$ its pointwise stabilizer (for short, we often
write $G_{x_1 \ldots x_m}$ in the latter case).

For \mbox{$\D=(X,\B,I)$} a Steiner \mbox{$t$-design} with \mbox{$G
\leq \Aut (\D)$}, let $G_B$ denote the setwise stabilizer of a block
$B \in \B$, and for $x \in X$, we define $G_{xB}= G_x \cap G_B$.

Let $\N$ be the set of positive integers (in this article, $0 \notin
\N$). For integers $m$ and $n$, let $(m,n)$ denote the greatest
common divisor of $m$ and $n$, and we write $m \mid n$ if $m$
divides $n$.

For any $x \in \R$, let $\lfloor x \rfloor$ denote the greatest
positive integer which is at most $x$.

All other notation is standard.

\medskip

When considering a Steiner \mbox{$t$-design} $\D$ with $t=2$, it is
elementary that the point \mbox{$2$-transitivity} of \mbox{$G \leq
\Aut(\D)$} implies its flag-transitivity. However, for $t \geq 3$,
it can be deduced from a result of
R.~E.~Block~\cite[Thm.\,2]{Block1965} that the converse holds:

\begin{prop}{\em (cf.~\cite{Buek1968,Hu2005}).}\label{flag2trs}
Let $\D=(X,\B,I)$ be a Steiner \mbox{$t$-design} with $t \geq 3$. If
\mbox{$G \leq \aut(\D)$} acts flag-transitively on $\D$, then $G$
also acts point \mbox{$2$-transitively} on $\D$.
\end{prop}

For $t \geq 5$, the flag-transitivity of \mbox{$G \leq \Aut(\D)$}
has an even stronger implication due to the following assertion,
which follows from Block's theorem and a combinatorial result of
D.~K.~Ray-Chaudhuri and R.~M.~Wilson~\cite[Thm.\,1]{Ray-ChWil1975}.

\begin{prop}{\em (cf.~\cite{CamPrae1993}).}\label{flag3hom}
Let $\D=(X,\B,I)$ be a Steiner \mbox{$t$-design} with $t\geq 2$.
Then, the following holds:

\begin{enumerate}

\item[(a)] If \mbox{$G \leq \aut(\D)$} acts block-transitively on $\D$,
then $G$ also acts point \linebreak \mbox{$\lfloor t/2
\rfloor$-homogeneously} on $\D$.

\item[(b)] If \mbox{$G \leq \aut(\D)$} acts flag-transitively on $\D$,
then $G$ also acts point \linebreak \mbox{$\lfloor (t+1)/2
\rfloor$-homogeneously} on $\D$.

\end{enumerate}
\end{prop}

We note that Propositions~\ref{flag2trs} and~\ref{flag3hom} hold
also for arbitrary $\lambda$, whereas for a
\mbox{$2$-$(v,k,\lambda)$} design the implication that the point
\mbox{$2$-transitivity} yields its flag-transitivity is only true if
\mbox{$(r,\lambda)=1$} (see,
e.g.,~\cite[Ch.\,2.3,\,Lemma\,8]{Demb1968}).

In order to investigate in the following all flag-transitive Steiner
\mbox{$5$-designs} and \mbox{$6$-designs}, we can as a consequence
of Proposition~\ref{flag3hom}~(b) make use of the classification of
all finite \mbox{$3$-homogeneous} permutation groups, which itself
relies on the classification of all finite simple groups
(cf.~\cite{Cam1981,Gor1982,Kant1972,Lieb1987,LivWag1965}).

The list of groups is as follows.

Let $G$ be a finite \mbox{$3$-homogeneous} permutation group on a
set $X$ with $\left|X\right| \geq 4$. Then $G$ is either of

{\bf (A) Affine Type:} $G$ contains a regular normal subgroup $T$
which is elementary Abelian of order $v=2^d$. If we identify $G$
with a group of affine transformations
\[x \mapsto x^g+u\]
of $V=V(d,2)$, where $g \in G_0$ and $u \in V$, then particularly
one of the following occurs:

\begin{enumerate}

\smallskip

\item[(1)] $G \cong AGL(1,8)$, $A \mathit{\Gamma} L(1,8)$, or $A \mathit{\Gamma} L(1,32)$

\smallskip

\item[(2)] $G_0 \cong SL(d,2)$, $d \geq 2$

\smallskip

\item[(3)] $G_0 \cong A_7$, $v=2^4$

\end{enumerate}

\smallskip

or

\medskip

{\bf (B) Almost Simple Type:} $G$ contains a simple normal subgroup
$N$, and \mbox{$N \leq G \leq \Aut(N)$}. In particular, one of the
following holds, where $N$ and $v=|X|$ are given as follows:
\begin{enumerate}

\smallskip

\item[(1)] $A_v$, $v \geq 5$

\smallskip

\item[(2)] $PSL(2,q)$, $q>3$, $v=q+1$

\smallskip

\item[(3)] $M_v$, $v=11,12,22,23,24$ \hfill (Mathieu groups)

\smallskip

\item[(4)] $M_{11}$, $v=12$

\end{enumerate}

\medskip

We note that if $q$ is odd, then $PSL(2,q)$ is $3$-homogeneous for
\mbox{$q \equiv 3$ (mod $4$)}, but not for \mbox{$q \equiv 1$ (mod
$4$)}, and hence not every group $G$ of almost simple type
satisfying (2) is $3$-homogeneous on $X$. For required basic
properties of the listed groups, we refer, e.g.,
to~\cite{Atlas1985},~\cite{HupI1967},~\cite[Ch.\,2,\,5]{KlLi1990}.

We will now recall some standard combinatorial results on which we
rely in the sequel. Let $r$ (respectively $\lambda_2$) denote the
total number of blocks incident with a given point (respectively
pair of distinct points), and let all further parameters be as
defined at the beginning of Section~\ref{intro}.

\begin{lem}\label{divprop}
Let $\D=(X,\B,I)$ be a Steiner \mbox{$t$-design}. If $G \leq
\aut(\D)$ acts flag-transitively on $\D$, then
\[r  \bigm|  \left| G_x \right|\]
for any $x \in X$.
\end{lem}

\begin{lem} \label{Comb_t=5}
Let $\D=(X,\B,I)$ be a \mbox{$t$-$(v,k,\lambda)$} design. Then the
following holds:
\begin{enumerate}

\smallskip

\item[(a)] $bk = vr.$

\smallskip

\item[(b)] $\displaystyle{{v \choose t} \lambda = b {k \choose
t}.}$

\smallskip

\item[(c)] $r(k-1)=\lambda_2(v-1)$ for $t \geq 2$, where
$\displaystyle{\lambda_2=\lambda \frac{{v-2 \choose t-2}}{{k-2
\choose t-2}}.}$

\end{enumerate}
\end{lem}

\medskip

\begin{prop}{\em (cf.~\cite{Cam1976,Tits1964}).}\label{Cam}
If $\D=(X,\B,I)$ is a non-trivial Steiner \mbox{$t$-design}, then
the following holds:
\begin{enumerate}

\smallskip

\item[(a)] $v\geq (t+1)(k-t+1).$

\smallskip

\item[(b)] $v-t+1 \geq (k-t+2)(k-t+1)$ for $t>2$. If equality
holds, then
\smallskip
$(t,k,v)=(3,4,8),(3,6,22),(3,12,112),(4,7,23)$, or $(5,8,24)$.
\end{enumerate}
\end{prop}

We note that (a) is stronger for $k<2(t-1)$, while (b) is stronger
for $k>2(t-1)$. For $k=2(t-1)$ both assert that $v \geq t^2-1$.

As we are in particular interested in the case when $3 \leq t \leq
6$, we deduce from (b) the following upper bound for the positive
integer $k$.

\begin{Cor}\label{Cameron_t=5}
Let $\D=(X,\B,I)$ be a non-trivial Steiner \mbox{$t$-design} with
\mbox{$t=3+i$}, where \mbox{$i=0,1,2,3$}. Then
\[k \leq \bigl\lfloor \sqrt{v} + \textstyle{\frac{3}{2}+i} \bigr\rfloor.\]
\end{Cor}

\smallskip

\begin{rem} \label{equa_t=5}
\emph{If \mbox{$G \leq \Aut(\D)$} acts flag-transitively on any
Steiner \mbox{$t$-design} $\D$ with $t\geq 3$, then applying
Proposition~\ref{flag2trs} and Lemma~\ref{Comb_t=5}~(b) yields the
equation
\[b=\frac{{v \choose t}}{{k \choose t}}=\frac{v(v-1)
\left|G_{xy}\right|}{\left| G_B \right|},\] where $x$ and $y$ are
two distinct points in $X$ and $B$ is a block in $\B$, and thus
\[{v-2 \choose t-2} = (k-1) {k-2 \choose t-2}
\frac{\left|G_{xy} \right|}{\left|G_{xB}\right|} \;\, \mbox{if} \;\,
x \in B.\]}
\end{rem}

\bigskip


\section{The Classification of Flag-transitive Steiner 5-Designs}\label{flag5designs}

The classification of all non-trivial Steiner \mbox{$5$-designs}
admitting a flag-transitive group of automorphisms is as follows.

\begin{mthm}\label{class5-des}
Let $\D=(X,\B,I)$ be a non-trivial Steiner \mbox{$5$-design}. Then
\mbox{$G \leq \aut(\D)$} acts flag-transitively on $\D$ if and only
if one of the following occurs:

\begin{enumerate}
\item[(1)] $\D$ is isomorphic to the Witt \mbox{$5$-$(12,6,1)$}
design, and \mbox{$G \cong M_{12}$},

\medskip

\item[(2)] $\D$ is isomorphic to the Witt \mbox{$5$-$(24,8,1)$}
design, and \mbox{$G \cong PSL(2,23)$} or \mbox{$G \cong M_{24}$}.
\end{enumerate}
\end{mthm}

\medskip

We remark that in Part (2), \mbox{$G \cong PSL(2,23)$} acts sharply
flag-transitively on $\D$, and furthermore that $M_{24}$ as the full
group of automorphisms of $\D$ contains only one conjugacy class of
subgroups isomorphic to $PSL(2,23)$ (cf.~\cite{Atlas1985}).

\bigskip


\subsection{Groups of Automorphisms of Affine
Type}\label{affine typ} \hfill

\bigskip

In this subsection, we start with the proof of Main
Theorem~\ref{class5-des}. Using the notation as before, \emph{let
$\D=(X,\B,I)$ be a non-trivial Steiner \mbox{$5$-design} with
\mbox{$G \leq \aut(\D)$} acting flag-transitively on $\D$ throughout
the proof}. We recall that due to Proposition~\ref{flag3hom}, we may
restrict ourselves to the consideration of the finite
\mbox{$3$-homogeneous} permutation groups listed in
Section~\ref{Prelim}. Clearly, in the following we may assume that
$k>5$ as trivial Steiner \mbox{$5$-designs} are excluded. Let us
first assume that $G$ is of affine type.


\bigskip
\emph{Case} (1): $G \cong AGL(1,8)$, $A \mathit{\Gamma}
L(1,8)$, or $A \mathit{\Gamma} L(1,32)$.
\medskip

We may assume that $k > 5$. For $v=8$, we obtain $k=6$ by
Corollary~\ref{Cameron_t=5}, which is not possible in view of
Lemma~\ref{Comb_t=5}~(b). If $v=32$, then $\left| G \right|
=5v(v-1)$, and Lemma~\ref{divprop} immediately yields that \mbox{$G
\leq \Aut(\D)$} cannot act flag-transitively on any non-trivial
Steiner \mbox{$5$-design} $\D$.

\bigskip
\emph{Case} (2): $G_0 \cong SL(d,2)$, $d \geq 2$.
\medskip

Let $e_i$ denote the $i$-th standard basis vector of the vector
space $V=V(d,2)$, and $\text{\footnotesize{$\langle$}} e_i
\text{\footnotesize{$\rangle$}}$ the \mbox{$1$-dimensional} vector
subspace spanned by $e_i$. We will prove by contradiction that
\mbox{$G \leq \Aut(\D)$} cannot act flag-transitively on any
non-trivial Steiner \mbox{$5$-design} $\D$.

We may assume that $v=2^d > k > 5$. For $d=3$, we have $v=8$ and
$k=6$ by Corollary~\ref{Cameron_t=5}, which is not possible in view
of Lemma~\ref{Comb_t=5}~(b) again. So, we may assume that $d>3$. We
remark that clearly any five distinct points are non-coplanar in
$AG(d,2)$ and hence generate an affine subspace of dimension at
least $3$. Let $\mathcal{E}= \text{\footnotesize{$\langle$}}
e_1,e_2,e_3 \text{\footnotesize{$\rangle$}}$ denote the
\mbox{$3$-dimensional} vector subspace spanned by $e_1,e_2,e_3$.
Then by linear algebra $SL(d,2)_\mathcal{E}$, and therefore also
$G_{0,\mathcal{E}}$, acts point-transitively on \mbox{$V \setminus
\mathcal{E}$}. If the unique block $B \in \B$ which is incident with
the \mbox{$5$-subset} $\{0,e_1,e_2,e_3,e_1+e_2\}$ contains some
point outside $\mathcal{E}$, then it would already contain all
points of \mbox{$V \setminus \mathcal{E}$}. But then, we would have
\mbox{$k \geq 2^d-8+5=2^d-3$}, a contradiction to
Corollary~\ref{Cameron_t=5}. Hence, $B$ lies completely in
$\mathcal{E}$, and by the flag-transitivity of $G$, it follows that
each block must be contained in a \mbox{$3$-dimensional} affine
subspace. Thus, clearly $k \leq 8$. But, on the other hand, for $\D$
to be a block-transitive \mbox{$5$-design} admitting \mbox{$G \leq
\Aut(\D)$}, we obtain from~\cite{Alltop1971} the necessary (and
sufficient) condition that $2^d-3$ must divide $k \choose 4$, and
hence it follows for each respective value of $k$ that $d=3$,
contradicting our assumption.

\bigskip
\emph{Case} (3): $G_0 \cong A_7$, $v=2^4$.
\medskip

Since $v=2^4$, we obtain from Corollary~\ref{Cameron_t=5} that $k
\leq 7$. But, Lemma~\ref{divprop} easily rules out the cases when
$k=6$ or $7$.

\bigskip


\subsection{\mbox{Groups of Automorphisms of Almost Simple
Type}}\label{almost simple type} \hfill

\bigskip

Before we consider in this subsection successively those cases where
$G$ is of almost simple type, we indicate some lemmas which will be
required for Case (2).


\smallskip

Let $q$ be a prime power $p^e$, and $U$ a subgroup of $PSL(2,q)$.
Furthermore, let $N_l$ denote the number of orbits of length $l$ and
let \linebreak \mbox{$n=(2,q-1)$}. In~\cite[Ch.\,5]{Hu_Habil2005},
we have in particular determined the orbit-lengths from the action
of subgroups of $PSL(2,q)$ on the points of the projective line. For
the list of subgroups of $PSL(2,q)$, we thereby refer
to~\cite[Ch.\,12,\,p.\,285f.]{Dick1901}
or~\cite[Ch.\,II,\,Thm.\,8.27]{HupI1967}.

\medskip

\begin{lem}\label{PSL_cyc}
Let $U$ be the cyclic group of order $c$ with $c \mid \frac{q \pm
1}{n}$. Then
\begin{enumerate}
\item[(a)] if $c \mid \frac{q+1}{n}$, then $N_c=(q+1)/c$,
\item[(b)] if $c \mid \frac{q-1}{n}$, then $N_1=2$ and $N_c=(q-1)/c$.
\end{enumerate}
\end{lem}

\begin{lem}\label{PSL_dihed}
Let $U$ be the dihedral group of order $2c$ with $c \mid \frac{q \pm
1}{n}$. Then
\begin{enumerate}
\item[(i)] for \mbox{$q \equiv 1$ $($\emph{mod} $4)$}:
\begin{enumerate}
\item[(a)] if $c \mid \frac{q+1}{2}$, then $N_c=2$ and $N_{2c}=(q+1-2c)/(2c)$,
\item[(b)] if $c \mid \frac{q-1}{2}$, then $N_2=1$, $N_c=2$, and $N_{2c}=(q-1-2c)/(2c)$,
unless $c=2$, in which case $N_2=3$ and $N_4=(q-5)/4$,
\end{enumerate}
\item[(ii)] for \mbox{$q \equiv 3$ $($\emph{mod} $4)$}:
\begin{enumerate}
\item[(a)] if $c \mid \frac{q+1}{2}$, then $N_{2c}=(q+1)/(2c)$,
\item[(b)] if $c \mid \frac{q-1}{2}$, then $N_2=1$ and
$N_{2c}=(q-1)/(2c)$,
\end{enumerate}
\item[(iii)] for \mbox{$q \equiv 0$ $($\emph{mod} $2)$}:
\begin{enumerate}
\item[(a)] if $c \mid q+1$, then $N_c=1$ and $N_{2c}=(q+1-c)/(2c)$,
\item[(b)] if $c \mid q-1$, then $N_2=1$, $N_c=1$, and $N_{2c}=(q-1-c)/(2c)$.
\end{enumerate}
\end{enumerate}
\end{lem}

\begin{lem}\label{PSL_elAb}
Let $U$ be the elementary Abelian group of order $\q \mid q$. Then
$N_1=1$ and $N_{\q}=q/ \q$.
\end{lem}

\begin{lem}\label{PSL_semi}
Let $U$ be a semi-direct product of the elementary Abelian group of
order $\q \mid q$ and the cyclic group of order $c$ with $c \mid
\q-1$ and $c \mid q-1$. Then $N_1=1$, $N_{\q}=1$, and $N_{c \q}=(q-
\q)/(c \q)$.
\end{lem}

\begin{lem}\label{PSL_kl.PSL}
Let $U$ be $PSL(2,\q)$ with $\q^m = q$, $m \geq 1$. Then
$N_{\q+1}=1$, $N_{\q (\q-1)}=1$ if $m$ is even, and all other orbits
are regular.
\end{lem}

\begin{lem}\label{PSL_kl.PGL}
Let $U$ be $PGL(2,\q)$ with $\q^m = q$, $m > 1$ even. Then
$N_{\q+1}=1$, $N_{\q(\q-1)}=1$, and all other orbits are regular.
\end{lem}

\begin{lem}\label{PSL_A_4}
Let $U$ be isomorphic to $A_4$. Then
\begin{enumerate}
\item[(i)] for \mbox{$q \equiv 1$ $($\emph{mod} $4)$}:
\begin{enumerate}
\item[(a)] if $3 \mid \frac{q+1}{2}$, then $N_{6}=1$ and $N_{12}=(q-5)/12$,
\item[(b)] if $3 \mid \frac{q-1}{2}$, then $N_4=2$, $N_{6}=1$, and $N_{12}=(q-13)/12$,
\item[(c)] if $3 \mid q$, then $N_4=1$, $N_{6}=1$, and $N_{12}=(q-9)/12$,
\end{enumerate}
\item[(ii)] for \mbox{$q \equiv 3$ $($\emph{mod} $4)$}:
\begin{enumerate}
\item[(a)] if $3 \mid \frac{q+1}{2}$, then $N_{12}=(q+1)/12$,
\item[(b)] if $3 \mid \frac{q-1}{2}$, then $N_4=2$ and $N_{12}=(q-7)/12$,
\item[(c)] if $3 \mid q$, then $N_4=1$ and $N_{12}=(q-3)/12$,
\end{enumerate}
\item[(iii)] for $q=2^e$, \mbox{$e \equiv 0$ $($\emph{mod} $2)$}:
$N_{1}=1$, $N_{4}=1$, and \\$N_{12}=(q-4)/12$.
\end{enumerate}
\end{lem}

\begin{lem}\label{PSL_S_4}
Let $U$ be isomorphic to $S_4$. Then
\begin{enumerate}
\item[(i)] for \mbox{$q \equiv 1$ $($\emph{mod} $8)$}:
\begin{enumerate}
\item[(a)] if $3 \mid \frac{q+1}{2}$, then $N_{6}=1$, $N_{12}=1$, and $N_{24}=(q-17)/24$,
\item[(b)] if $3 \mid \frac{q-1}{2}$, then $N_{6}=1$, $N_{8}=1$, $N_{12}=1$, and \\$N_{24}=(q-25)/24$,
\item[(c)] if $3 \mid q$, then $N_4=1$, $N_{6}=1$, and $N_{24}=(q-9)/24$,
\end{enumerate}
\item[(ii)] for \mbox{$q \equiv -1$ $($\emph{mod} $8)$}:
\begin{enumerate}
\item[(a)] if $3 \mid \frac{q+1}{2}$, then $N_{24}=(q+1)/24$,
\item[(b)] if $3 \mid \frac{q-1}{2}$, then $N_8=1$ and
$N_{24}=(q-7)/24$.
\end{enumerate}
\end{enumerate}
\end{lem}

\begin{lem}\label{PSL_A_5}
Let $U$ be isomorphic to $A_5$. Then
\begin{enumerate}
\item[(i)] for \mbox{$q \equiv 1$ $($\emph{mod} $4)$}:
\begin{enumerate}
\item[(a)] if $q=5^e$, \mbox{$e \equiv 1$ $($\emph{mod} $2)$}, then $N_{6}=1$ and $N_{60}=(q-5)/60$,
\item[(b)] if $q=5^e$, \mbox{$e \equiv 0$ $($\emph{mod} $2)$}, then $N_{6}=1$, $N_{20}=1$, and \\ $N_{60}=(q-25)/60$,
\item[(c)] if $15 \mid \frac{q+1}{2}$, then $N_{30}=1$ and $N_{60}=(q-29)/60$,
\item[(d)] if $3 \mid \frac{q+1}{2}$ and $5 \mid \frac{q-1}{2}$, then $N_{12}=1$, $N_{30}=1$, and \\ $N_{60}=(q-41)/60$,
\item[(e)] if $3 \mid \frac{q-1}{2}$ and $5 \mid \frac{q+1}{2}$, then $N_{20}=1$, $N_{30}=1$, and \\ $N_{60}=(q-49)/60$,
\item[(f)] if $15 \mid \frac{q-1}{2}$, then $N_{12}=1$, $N_{20}=1$, $N_{30}=1$, and \\ $N_{60}=(q-61)/60$,
\item[(g)] if $3 \mid q$ and $5 \mid \frac{q+1}{2}$, then $N_{10}=1$ and $N_{60}=(q-9)/60$,
\item[(h)] if $3 \mid q$ and $5 \mid \frac{q-1}{2}$, then $N_{10}=1$, $N_{12}=1$, and \\$N_{60}=(q-21)/60$,
\end{enumerate}
\item[(ii)] for \mbox{$q \equiv 3$ $($\emph{mod} $4)$}:
\begin{enumerate}
\item[(a)] if $15 \mid \frac{q+1}{2}$, then $N_{60}=(q+1)/60$,
\item[(b)] if $3 \mid \frac{q+1}{2}$ and $5 \mid \frac{q-1}{2}$, then $N_{12}=1$ and $N_{60}=(q-11)/60$,
\item[(c)] if $3 \mid \frac{q-1}{2}$ and $5 \mid \frac{q+1}{2}$, then $N_{20}=1$ and $N_{60}=(q-19)/60$,
\item[(d)] if $15 \mid \frac{q-1}{2}$, then $N_{12}=1$, $N_{20}=1$, and
$N_{60}=(q-31)/60$.
\end{enumerate}
\end{enumerate}
\end{lem}

\medskip

We shall now turn to the examination of those cases where \mbox{$G
\leq \Aut(\D)$} is of almost simple type.

\bigskip
\emph{Case} (1): $N=A_v$, $v \geq 5$.
\medskip

We may assume that $v \geq 7 $. But then $A_v$, and hence also $G$,
is \mbox{$5$-transitive} and does not act on any non-trivial Steiner
\mbox{$5$-design} $\D$ in view of~\cite[Thm.\,3]{Kant1985}.

\bigskip
\emph{Case} (2): $N=PSL(2,q)$, $v=q+1$, $q=p^e >3$.
\medskip

Here $\Aut(N)= P \mathit{\Gamma} L (2,q)$, and $\left| G \right| =
(q+1)q \frac{(q-1)}{n}a$ with $n=(2,q-1)$ and $a \mid ne$. We may
assume that $q \geq 5$. We will show that only the flag-transitive
design given in Part (2) of Main Theorem~\ref{class5-des} with
\mbox{$G \cong PSL(2,23)$} can occur.

\emph{We will first assume that $N=G$.} Then, by
Remark~\ref{equa_t=5}, we obtain
\begin{equation}\label{Eq-0}
(q-2)(q-3) \left| PSL (2,q)_{0B} \right| \cdot n  =
(k-1)(k-2)(k-3)(k-4).
\end{equation}
In view of Proposition~\ref{Cam}~(b), we have
\begin{equation}\label{Eq-A}
q-3 \geq (k-3)(k-4),
\end{equation}
and thus it follows from equation~(\ref{Eq-0}) that
\begin{equation}\label{Eq-B}
(q-2)\left| PSL (2,q)_{0B} \right| \cdot n \leq (k-1)(k-2).
\end{equation}
If we assume that $k \geq 9$, then clearly
\[(k-1)(k-2) < 2(k-3)(k-4),\]
and hence we obtain
\[(q-2)\left| PSL (2,q)_{0B} \right| \cdot n< 2(q-3)\]
due to Proposition~\ref{Cam}~(b) again, which is obviously only
possible when\linebreak $\left| PSL (2,q)_{0B} \right| \cdot n=1$.
Thus, in particular $q$ has to be even. But then, considering
equation~(\ref{Eq-0}) yields that the left hand side of the equation
is not divisible by $4$, whereas obviously the right hand side is
always divisible by $8$, a contradiction. If $k<9$, then, using
equation~(\ref{Eq-0}) and inequality~(\ref{Eq-A}), the very few
remaining possibilities for $k$ can immediately be ruled out by
hand, except for the case when $k=8$, $q=23$ and $\left| PSL
(2,q)_{0B} \right|=1$. It is well-known that for the parameters
$t=5$, $v=24$ and $k=8$ there exists (up to isomorphism) only the
unique Witt \mbox{$5$-$(24,8,1)$} design $\D$, which can be
constructed from $PSL(2,23)$ in its natural \mbox{$3$-homogeneous}
action on the elements of \mbox{$GF(23) \cup \{\infty\}$}.
Furthermore, it can be shown that the setwise stabilizer
$PSL(2,23)_B$ of an appropriate, unique block $B \in \B$ is a
dihedral group of order $8$ (see,
e.g.,~\cite[Ch.\,IV,\,1.5]{BJL1999},~\cite[Ch.\,XIV,\,115]{Carm1937},
and~\cite[Thm.\,5]{Witt1938} for a uniqueness proof). Thus, using
Lemma~\ref{Comb_t=5}~(b), we obtain $b=759=\big[PSL(2,23) :
PSL(2,23)_{B}\big]$, and hence $PSL(2,23)$ acts block-transitively
on $\D$. As for $q=23$, the dihedral group of order $8$ has only
orbits of length $8$ in view of Lemma~\ref{PSL_dihed}~(ii)(a),
clearly $PSL(2,23)_B$ acts transitively on the points of $B$. Since
we have $\left| PSL (2,23)_{0B} \right|=1$, it follows that
$PSL(2,23)$ acts even sharply flag-transitively on $\D$.

\emph{Now, let us assume that $N<G \leq \aut(N)$.} We recall that
$q=p^e>3$, and will distinguish in the following the cases $p>3$,
$p=2$, and $p=3$.

\emph{First, let $p>3$.} We define $G^*=G \cap (PSL(2,q) \rtimes
\text{\footnotesize{$\langle$}} \tau_\alpha
\text{\footnotesize{$\rangle$}})$ with $\tau_\alpha \in$
Sym$(GF(p^e) \cup \{\infty\}) \cong S_v$ of order $e$ induced by the
Frobenius automorphism $\alpha : GF(p^e) \longrightarrow GF(p^e),\,
x \mapsto x^p$. Then, by Dedekind's law, we can write
\begin{equation}\label{Eq_G^*0}
G^*= PSL(2,q) \rtimes (G^* \cap \text{\footnotesize{$\langle$}}
\tau_\alpha \text{\footnotesize{$\rangle$}}).
\end{equation}
Defining $P \Si L(2,q)= PSL(2,q) \rtimes
\text{\footnotesize{$\langle$}} \tau_\alpha
\text{\footnotesize{$\rangle$}}$, it can easily be calculated that
$P \Si L(2,q)_{0,1,\infty} = \text{\footnotesize{$\langle$}}
\tau_\alpha \text{\footnotesize{$\rangle$}}$, and
$\text{\footnotesize{$\langle$}} \tau_\alpha
\text{\footnotesize{$\rangle$}}$ has precisely $p+1$ distinct fixed
points (cf., e.g.,~\cite[Ch.\,6.4,\,Lemma\,2]{Demb1968}). As $p>3$,
we conclude therefore that $G^* \cap \text{\footnotesize{$\langle$}}
\tau_\alpha \text{\footnotesize{$\rangle$}} \leq G^*_{0B}$ for some
appropriate, unique block $B \in \B$ by the definition of Steiner
\mbox{$5$-designs}. Furthermore, clearly $PSL(2,q) \cap (G^* \cap
\text{\footnotesize{$\langle$}} \tau_\alpha
\text{\footnotesize{$\rangle$}}) =1.$ Hence, we have
\begin{equation}\label{Eq_G^*}
\begin{split}
\left| (0,B)^{G^*} \right| &=  \big[G^* : G^*_{0B} \big] \\
                           &= \big[PSL(2,q) \rtimes (G^* \cap \text{\footnotesize{$\langle$}}
\tau_\alpha \text{\footnotesize{$\rangle$}}) : PSL(2,q)_{0B} \rtimes
(G^* \cap \text{\footnotesize{$\langle$}} \tau_\alpha
\text{\footnotesize{$\rangle$}})\big] \\
                           &= \big[PSL(2,q) : PSL(2,q)_{0B} \big] \\
                           & =  \left| (0,B)^{PSL(2,q)}\right|.\\
\end{split}
\end{equation}
Thus, if we assume that \mbox{$G^* \leq \Aut(\D)$} acts already
flag-transitively on $\D$, then we obtain $\left| (0,B)^{G^*}
\right|=\left| (0,B)^{PSL(2,q)}\right|=bk$ in view of
Remark~\ref{equa_t=5}. Hence, $PSL(2,q)$ must also act
flag-transitively on $\D$, and we may proceed as in the case when
$N=G$. Therefore, let us assume that \mbox{$G^* \leq \Aut(\D)$} does
not act flag-transitively on $\D$. Then, we conclude that $\big[G :
G^* \big]=2$ and $G^*$ has exactly two orbits of equal length on the
set of flags. Thus, by equation~(\ref{Eq_G^*}), we obtain for the
orbit containing the flag $(0,B)$ that $\left| (0,B)^{G^*}
\right|=\left| (0,B)^{PSL(2,q)}\right|=\frac{bk}{2}$. As it is
well-known the normalizer of $PSL(2,q)$ in Sym$(X)$ is $P \G L
(2,q)$, and hence in particular $PSL(2,q)$ is normal in $G$. It
follows therefore that we have under $PSL(2,q)$ also precisely one
further orbit of equal length on the set of flags. Then, proceeding
similarly to the case $N=G$ for each orbit on the set of flags, we
obtain (representative for the orbit containing the flag $(0,B)$)
that
\begin{equation}\label{Eq-0_N<G}
\frac{(q-2)(q-3) \left| PSL (2,q)_{0B} \right| \cdot n}{2}  =
(k-1)(k-2)(k-3)(k-4),
\end{equation}
and as here $n=2$, this is is equivalent to
\begin{equation}\label{Eq-0_N<G-equiv}
\begin{split}
(q-2)(q-3) \left| PSL (2,q)_{0B} \right| & =  (k-1)(k-2)(k-3)(k-4)\\
& =  k(k^3-10k^2+35k-50)+24.
\end{split}
\end{equation}
Hence, we have in particular
\begin{equation}\label{Eq-1_N<G}
k \bigm| (q-2)(q-3) \left| PSL (2,q)_{0B} \right| -24.
\end{equation}
Since $PSL(2,q)_B$ can have one or two orbits of equal length on the
set of points of $B$, we have
\begin{equation}\label{Eq-k_N<G}
k\;\,\mbox{or}\;\,\frac{k}{2}=\left|0^{PSL(2,q)_B}\right| = \big[PSL(2,q)_B : PSL(2,q)_{0B}
\big].
\end{equation}
By the same arguments as in case $N=G$, we deduce from
equation~(\ref{Eq-0_N<G-equiv}) that
\begin{equation}\label{Eq-C}
(q-2)\left| PSL (2,q)_{0B} \right| \leq (k-1)(k-2),
\end{equation}
and assuming that $k \geq 9$, we obtain
\[(q-2)\left| PSL (2,q)_{0B} \right| < 2(q-3),\]
which is clearly only possible when \mbox{$\left| PSL (2,q)_{0B}
\right|=1$.} Hence, it follows that
\begin{equation}\label{Eq-D}
(q-2)(q-3)=(k-1)(k-2)(k-3)(k-4),
\end{equation}
and $k \mid(q-2)(q-3)-24$ in view of property~(\ref{Eq-1_N<G}). On
the other hand, for \mbox{$k \geq 9$}, we obtain from
equation~(\ref{Eq-k_N<G}) that $k$ or
$\frac{k}{2}=\left|{PSL(2,q)_B}\right| \bigm|
\left|PSL(2,q)\right|=\frac{q^3-q}{2}$, and thus in particular $k
\mid q^3-q$. But, it can easily be seen that $(q^3-q,(q-2)(q-3)-24)
\mid 2^3\cdot 3 \cdot 11$, and thus we have only a small number of
possibilities for $k$ to check, which can easily be eliminated by
hand using equation~(\ref{Eq-D}). For $k<9$, the very few remaining
possibilities for $k$ can immediately be ruled out by hand using
inequality~(\ref{Eq-A}) and equation~(\ref{Eq-0_N<G-equiv}), except
for the case when $k=8$, $q=23$ and $\left| PSL (2,q)_{0B}
\right|=2$. But, as involutions are fixed point free on the points
of the projective line for \mbox{$q \equiv 3$ (mod $4$)} (cf.,
e.g.,~\cite[Ch.\,II,\,Thm.\,8.5]{HupI1967}), this is impossible.

\emph{Now, let $p=2$.} Then, clearly $N=PSL(2,q)=PGL(2,q)$, and we
have $\Aut(N)=P\Si L(2,q)$. If we assume that
$\text{\footnotesize{$\langle$}} \tau_\alpha
\text{\footnotesize{$\rangle$}} \leq P \Si L (2,q)_{0B}$ for some
appropriate, unique block $B \in \B$, then, using the terminology
of~(\ref{Eq_G^*0}), we have $G^*=G=P\Si L(2,q)$ and as clearly
$PSL(2,q) \cap \text{\footnotesize{$\langle$}} \tau_\alpha
\text{\footnotesize{$\rangle$}}=1$, we can apply
equation~(\ref{Eq_G^*}). Thus, $PSL(2,q)$ must also be
flag-transitive, which has already been considered. Therefore, we
may assume that $\text{\footnotesize{$\langle$}} \tau_\alpha
\text{\footnotesize{$\rangle$}} \nleq P\Si L(2,q)_{0B}$. Let $s$ be
a prime divisor of $e=\left| \text{\footnotesize{$\langle$}}
\tau_\alpha \text{\footnotesize{$\rangle$}} \right|$. As the normal
subgroup $H:=(P \Si L (2,q)_{0,1,\infty})^s \leq
\text{\footnotesize{$\langle$}} \tau_\alpha
\text{\footnotesize{$\rangle$}}$ of index $s$ has precisely $p^s+1$
distinct fixed points (see,
e.g.,~\cite[Ch.\,6.4,\,Lemma\,2]{Demb1968}), we have $G \cap H \leq
G_{0B}$ for some appropriate, unique block $B \in \B$ by the
definition of Steiner \mbox{$5$-designs}. It can then be deduced
that $e=s^u$ for some $u \in \N$, since if we assume for $G= P \Si
L(2,q)$ that there exists a further prime divisor $\overline{s}$ of
$e$ with $\overline{s} \neq s$, then $\overline{H}:=(P \Si L
(2,q)_{0,1,\infty})^{\overline{s}} \leq
\text{\footnotesize{$\langle$}} \tau_\alpha
\text{\footnotesize{$\rangle$}}$ and $H$ are both subgroups of $P\Si
L(2,q)_{0B}$ by the flag-transitivity of $P \Si L (2,q)$, and hence
$\text{\footnotesize{$\langle$}} \tau_\alpha
\text{\footnotesize{$\rangle$}} \leq P\Si L(2,q)_{0B}$, a
contradiction. Furthermore, as $\text{\footnotesize{$\langle$}}
\tau_\alpha \text{\footnotesize{$\rangle$}}\nleq P \Si L
(2,q)_{0B}$, we may, by applying Dedekind's law, assume that
\[G_{0B} = PSL(2,q) _{0B} \rtimes (G \cap H).\]
Thus, by Remark~\ref{equa_t=5}, we obtain
\[(q-2)(q-3) \left| PSL (2,q)_{0B} \right| \left| G \cap H \right| =
(k-1)(k-2)(k-3)(k-4) \left| G \cap \text{\footnotesize{$\langle$}}
\tau_\alpha \text{\footnotesize{$\rangle$}}\right|.\] Using that
$k=\left|0^{G_B}\right| = \big[G_B : G_{0B} \big]$, we have more
precisely
\begin{enumerate}
\item[(A)] if $G= PSL(2,q) \rtimes (G \cap H)$:
\[(q-2) (q-3)\left| PSL (2,q)_{0B} \right| = (k-1)(k-2)(k-3)(k-4)\]
\[\mbox{with} \;\, \left| PSL(2,q)_{0B}\right| =\frac{\left| PSL(2,q)_B \right|}{k}, \;\, \mbox{or}\]

\item[(B)] if $G = P \Si L (2,q)$:
\[(q-2)(q-3) \left| PSL (2,q)_{0B} \right| =(k-1)(k-2)(k-3)(k-4) s\]
\[\mbox{with} \;\, \left| PSL(2,q)_{0B}\right| =\frac{\left| PSL(2,q)_B
\right|}{k} \cdot \left\{\begin{array}{ll}
    s,\;\, \mbox{if} \;\, G_B=PSL (2,q)_B \rtimes \text{\footnotesize{$\langle$}} \tau_\alpha
\text{\footnotesize{$\rangle$}}\\
    1,\;\,\mbox{if} \;\, G_B=PSL (2,q)_B \rtimes H.\\
\end{array} \right.\]
\end{enumerate}

As far as condition~(A) is concerned, we may argue exactly as in the
earlier case $N=G$. Thus, only condition~(B) has to be examined, and
we will also show that here \mbox{$G \leq \Aut(\D)$} cannot act
flag-transitively on any non-trivial Steiner \mbox{$5$-design} $\D$.
Clearly, there exists always a Klein four-group $V_4 \leq PSL(2,q)$,
which fixes some \mbox{$4$-subset} $S$ of $X$ and some additional
point $x \in X$, and hence must fix the unique block $B \in \B$
which is incident with $S \cup \{x\}$ by the definition of Steiner
\mbox{$5$-designs}. Examining the list of possible subgroups of
$PSL(2,q)$ with their orbits on the projective line
(cf.~Lemmas~\ref{PSL_cyc}-\ref{PSL_A_5}), it follows that we only
have to consider the possibility when $PSL(2,q)_B$ is conjugate to
$PSL(2,\q)$ with $\q^m = q$, $m \geq 1$, and by
Lemma~\ref{PSL_kl.PSL}, we conclude that $k=\q+1$. Applying
condition~(B) yields then
\begin{equation}\label{eq-0}
(q-2)(q-3) \left| PSL(2,q)_{0B}\right| = \q(\q-1) (\q -2) (\q-3)s
\end{equation}
\[\mbox{with} \;\, \left| PSL(2,q)_{0B}\right| = \q (\q- 1)\cdot \left\{\begin{array}{ll}
    s,\;\mbox{or}\\
    1.\\
\end{array} \right.\]
Since $q=2^{s^u}$, we can write $\q=2^{s^w}$ for some integer $0
\leq w \leq u$, and $q=\q^m=\q^{s^{u-w}}$. As we may assume that
$k=\q +1 = 2^{s^w} +1>5$, it follows in particular that $w \geq 1$,
and hence $s < 2^{s^w}=\q$. Thus, using equation~(\ref{eq-0}), we
obtain
\[(\q^{s^{u-w}}-2)(\q^{s^{u-w}}-3)=(q-2)(q-3) \leq (\q - 2)(\q-3)s < ({\q}^2 -2s)(\q-3).\]
But, as clearly $u-w \geq 1$ (otherwise, $k=q+1$, a contradiction to
Corollary~\ref{Cameron_t=5}), this yields a contradiction for every
prime $s$.

\emph{Now, let $p=3$.} We have $\Aut(N)=P\Gamma L(2,q)=PGL(2,q)
\rtimes \text{\footnotesize{$\langle$}} \tau_\alpha
\text{\footnotesize{$\rangle$}}$, and as $PGL(2,q)$ is sharply
\mbox{$3$-transitive}, it follows that $P \Gamma L
(2,q)_{0,1,\infty} = \text{\footnotesize{$\langle$}} \tau_\alpha
\text{\footnotesize{$\rangle$}}$. Again, we define $G^*=G \cap
(PSL(2,q) \rtimes \text{\footnotesize{$\langle$}} \tau_\alpha
\text{\footnotesize{$\rangle$}})$ and may write $G^*= PSL(2,q)
\rtimes (G^* \cap \text{\footnotesize{$\langle$}} \tau_\alpha
\text{\footnotesize{$\rangle$}})$ as in equation~(\ref{Eq_G^*0}). We
distinguish the cases $G=G^*$ and \mbox{$[G:G^*]=2$}. In the
following, we will examine the first case in detail, whereas the
second may be treated mutatis mutandis. Let $G=G^*$. Then, we have
$\Aut(N)=P\Si L(2,q)$. If we assume that
$\text{\footnotesize{$\langle$}} \tau_\alpha
\text{\footnotesize{$\rangle$}} \leq P \Si L (2,q)_{0B}$ for some
appropriate, unique block $B \in \B$, then $G=P\Si L(2,q)$, and as
clearly $PSL(2,q) \cap \text{\footnotesize{$\langle$}} \tau_\alpha
\text{\footnotesize{$\rangle$}}=1$, we can apply
equation~(\ref{Eq_G^*}). Thus, $PSL(2,q)$ must also be
flag-transitive, which has already been considered. Therefore, we
may assume that $\text{\footnotesize{$\langle$}} \tau_\alpha
\text{\footnotesize{$\rangle$}} \nleq P\Si L(2,q)_{0B}$. Let $s$ be
a prime divisor of $e=\left| \text{\footnotesize{$\langle$}}
\tau_\alpha \text{\footnotesize{$\rangle$}} \right|$. As already
mentioned, the normal subgroup $H:=(P \Si L (2,q)_{0,1,\infty})^s
\leq \text{\footnotesize{$\langle$}} \tau_\alpha
\text{\footnotesize{$\rangle$}}$ of index $s$ has precisely $p^s+1$
distinct fixed points, and hence we have $G \cap H \leq G_{0B}$ for
some appropriate, unique block $B \in \B$ by the definition of
Steiner \mbox{$5$-designs}. It can then be deduced exactly as for
$p=2$ that $e=s^u$ for some $u \in \N$. As
$\text{\footnotesize{$\langle$}} \tau_\alpha
\text{\footnotesize{$\rangle$}}\nleq P \Si L (2,q)_{0B}$, we may, by
applying Dedekind's law, assume that
\[G_{0B} = PSL(2,q) _{0B} \rtimes (G \cap H).\]
Thus, by Remark~\ref{equa_t=5}, we obtain
\[2(q-2)(q-3) \left| PSL (2,q)_{0B} \right| \left| G \cap H \right|=
(k-1)(k-2)(k-3)(k-4) \left| G \cap \text{\footnotesize{$\langle$}}
\tau_\alpha \text{\footnotesize{$\rangle$}}\right|.\] Using that
$k=\left|0^{G_B}\right| = \big[G_B : G_{0B} \big]$, we have more
precisely
\begin{enumerate}
\item[(A$^*$)] if $G= PSL(2,q) \rtimes (G \cap H)$:
\[2(q-2)(q-3) \left| PSL (2,q)_{0B} \right| = (k-1)(k-2)(k-3)(k-4)\]
\[\mbox{with} \;\, \left| PSL(2,q)_{0B}\right| =\frac{\left| PSL(2,q)_B \right|}{k}, \;\, \mbox{or}\]

\item[(B$^*$)] if $G = P \Si L (2,q)$:
\[2(q-2)(q-3) \left| PSL (2,q)_{0B} \right| =(k-1)(k-2)(k-3)(k-4) s\]
\[\mbox{with} \;\, \left| PSL(2,q)_{0B}\right| =\frac{\left| PSL(2,q)_B
\right|}{k} \cdot \left\{\begin{array}{ll}
    s,\;\, \mbox{if} \;\, G_B=PSL (2,q)_B \rtimes \text{\footnotesize{$\langle$}} \tau_\alpha
\text{\footnotesize{$\rangle$}}\\
    1,\;\,\mbox{if} \;\, G_B=PSL (2,q)_B \rtimes H.\\
\end{array} \right.\]
\end{enumerate}

Considering condition~(A$^*$), we may argue exactly as in the
earlier case $N=G$. Thus, only condition~(B$^*$) has to be examined,
and we will show in the following that here \mbox{$G \leq \Aut(\D)$}
cannot act flag-transitively on any non-trivial Steiner
\mbox{$5$-design} $\D$. In view of the subgroups of $PSL(2,q)$ with
their orbits on the projective line
(Lemmas~\ref{PSL_cyc}-\ref{PSL_A_5}), we have to examine the
following possibilities:

\begin{enumerate}
\item[(i)] $PSL(2,q)_B$ is conjugate to a cyclic subgroup of order $c$
with $c \mid \frac{q \pm 1}{2}$ of $PSL(2,q)$, and $k=c$.

\item[(ii)] $PSL(2,q)_B$ is conjugate to a dihedral subgroup of order $2c$
with $c \mid \frac{q \pm 1}{2}$ of $PSL(2,q)$, and $k=c$ or $2c$.

\item[(iii)] $PSL(2,q)_B$ is conjugate to an elementary Abelian subgroup
of order $\q \mid q$ of $PSL(2,q)$, and $k=\q$.

\item[(iv)] $PSL(2,q)_B$ is conjugate to a semi-direct product
of an elementary Abelian subgroup of order $\q \mid q$ with a cyclic
subgroup of order $c$ of $PSL(2,q)$ with $c \mid \q-1$ and $c \mid
q-1$, and $k=\q$ or $c\q$.

\item[(v)] $PSL(2,q)_B$ is conjugate to $PSL(2,\q)$ with
$\q^m = q$, $m \geq 1$, and $k=\q+1$, $\q(\q-1)$ if $m$ is even, or
$k=(\q +1)\q (\q -1)/2$.

\item[(vi)] $PSL(2,q)_B$ is conjugate to $PGL(2,\q)$ with
$\q^m = q$, $m > 1$ even, and $k= \q +1$, $\q(\q-1)$ or $k=(\q +1)\q
(\q -1)$.

\item[(vii)] $PSL(2,q)_B$ is conjugate to $A_4$, and $k=6$ or $12$.

\item[(viii)] $PSL(2,q)_B$ is conjugate to $S_4$, and $k=6$ or $24$.

\item[(ix)] $PSL(2,q)_B$ is conjugate to $A_5$, and $k=10$, $12$ or $60$.
\end{enumerate}

Since $q=3^{s^u}$, we can write $\q=3^{s^w}$ for some integer $0
\leq w \leq u$, and $q=\q^m=\q^{s^{u-w}}$.

ad (i): By condition~(B$^*$), we have
\[2(q-2)(q-3) \left| PSL(2,q)_{0B}\right| = (c-1)(c-2)(c-3)(c-4)s\]
\[\mbox{with} \;\, \left| PSL(2,q)_{0B}\right| = \left\{\begin{array}{ll}
    s,\;\mbox{or}\\
    1.\\
\end{array} \right.\]
In view of the earlier case $N=G$, it is sufficient to consider the
equation
\begin{equation}\label{E1}
(q-2)(q-3)= \frac{(c-1)(c-2)(c-3)(c-4)s}{2}.
\end{equation}
For $c \mid \frac{q+1}{2}$, equation~(\ref{E1}) yields
\begin{equation*}
\begin{split}
c \bigm| \frac{(q+1)(q-6)}{2} &=\frac{(q-2)(q-3)}{2}-6 =
\frac{(c-1)(c-2)(c-3)(c-4)s}{4}-6\\ &=\frac{cs}{4}(c^3-10
c^2+35c-50) + 6s -6,\end{split}
\end{equation*}
and thus $c \mid 6s-6$ must hold. If $c \mid \frac{q-1}{2}$, then,
by equation~(\ref{E1}), we have
\begin{equation*}
\begin{split}
c \bigm| \frac{(q-1)(q-4)}{2} &=\frac{(q-2)(q-3)}{2}-1 =
\frac{(c-1)(c-2)(c-3)(c-4)s}{4}-1\\ &=\frac{cs}{4}(c^3-10
c^2+35c-50) + 6s -1,\end{split}
\end{equation*}
and hence $c \mid 6s-1$ must hold. As clearly $c<6s$ in both cases,
it follows from equation~(\ref{E1}) that in particular
\[(3^{s^u}-2)(3^{s^u-1}-1)<\frac{c^4 s}{6}<6^3 \cdot s^5,\]
which implies that $s^u \leq 7$. As $c \mid 6s-6$ respectively $c
\mid 6s-1$, this leaves only a very small number of possibilities
for $k$ to check, which can easily be ruled out by hand using
equation~(\ref{E1}).

ad (ii): Let $k=c$. Applying condition~(B$^*$) yields
\[2(q-2)(q-3) \left| PSL(2,q)_{0B}\right| = (c-1)(c-2)(c-3)(c-4)s\]
\[\mbox{with} \;\, \left| PSL(2,q)_{0B}\right| = 2 \cdot \left\{\begin{array}{ll}
    s,\;\mbox{or}\\
    1.\\
\end{array} \right.\]
First, let $k=c$. Due to the earlier case $N=G$, it is sufficient to
consider the equation
\begin{equation}\label{E2}
(q-2)(q-3)= \frac{(c-1)(c-2)(c-3)(c-4)s}{4}.
\end{equation}
If $c \mid \frac{q+1}{2}$, then, by equation~(\ref{E2}), we have
\begin{equation*}
\begin{split}
c \bigm| \frac{(q+1)(q-6)}{2} &=\frac{(q-2)(q-3)}{2}-6 =
\frac{(c-1)(c-2)(c-3)(c-4)s}{8}-6\\ &=\frac{cs}{8}(c^3-10
c^2+35c-50) + 3s -6,\end{split}
\end{equation*}
and hence $c \mid 3s-6$ must hold. For $c \mid \frac{q-1}{2}$, it
follows from equation~(\ref{E2}) that
\begin{equation*}
\begin{split}
c \bigm| \frac{(q-1)(q-4)}{2} &=\frac{(q-2)(q-3)}{2}-1 =
\frac{(c-1)(c-2)(c-3)(c-4)s}{8}-1\\ &=\frac{cs}{8}(c^3-10
c^2+35c-50) + 3s -1,\end{split}
\end{equation*}
and thus $c \mid 3s-1$ must hold. Obviously, we have $c<3s$ in both
cases, and therefore equation~(\ref{E2}) gives in particular
\[4(3^{s^u}-2)(3^{s^u-1}-1)<\frac{c^4 s}{3}<3^3 \cdot s^5,\]
which implies that $s^u \leq 5$. Due to the fact that $c \mid 3s-6$
respectively $c \mid 3s-1$, we have only a very small number of
possibilities for $k$ to check, which can easily be ruled out by
hand using equation~(\ref{E2}). Now, let $k=2c$. Due to
condition~(B$^*$), we have
\[2(q-2)(q-3) \left| PSL(2,q)_{0B}\right| = (2c-1)(2c-2)(2c-3)(2c-4)s\]
\[\mbox{with} \;\, \left| PSL(2,q)_{0B}\right| = \left\{\begin{array}{ll}
    s,\;\mbox{or}\\
    1.\\
\end{array} \right.\]
Again, it suffices to consider the equation
\begin{equation}\label{E3}
\frac{(q-2)(q-3)}{2}= (2c-1)(c-1)(2c-3)(c-2)s.
\end{equation}
For $c \mid \frac{q+1}{2}$, equation~(\ref{E3}) yields
\begin{equation*}
\begin{split}
c \bigm| \frac{(q+1)(q-6)}{2} &=\frac{(q-2)(q-3)}{2}-6 =
(2c-1)(c-1)(2c-3)(c-2)s-6\\ &=cs(4c^3-20 c^2+35c-25) + 6s
-6,\end{split}
\end{equation*}
and thus $c \mid 6s-6$ must hold. If $c \mid \frac{q-1}{2}$, then
due to equation~(\ref{E3}), we have
\begin{equation*}
\begin{split}
c \bigm| \frac{(q-1)(q-4)}{2} &=\frac{(q-2)(q-3)}{2}-1 =
(2c-1)(c-1)(2c-3)(c-2)s-1\\ &=cs(4c^3-20 c^2+35c-25) + 6s
-1,\end{split}
\end{equation*}
and hence $c \mid 6s-1$ must hold. As clearly $c<6s$ in both cases,
we deduce from equation~(\ref{E3}) that in particular
\[(3^{s^u}-2)(3^{s^u-1}-1)<\frac{(2c)^4 s}{6}<2^4 \cdot 6^3 \cdot s^5,\]
and hence it follows that $s^u \leq 7$. Since we have  $c \mid 6s-6$
respectively $c \mid 6s-1$, this leaves only a very small number of
possibilities for $k$ to check, which can easily be ruled out by
hand using equation~(\ref{E3}).

ad (iii): In view of condition~(B$^*$), we have
\[2(q-2)(q-3) \left| PSL(2,q)_{0B}\right| = (\q-1)(\q-2)(\q-3)(\q-4)s\]
\[\mbox{with} \;\, \left| PSL(2,q)_{0B}\right| = \left\{\begin{array}{ll}
    s,\;\mbox{or}\\
    1.\\
\end{array} \right.\]
It suffices to consider the equation
\begin{equation}\label{E4}
2(q-2)(q-3)= (\q-1)(\q-2)(\q-3)(\q-4)s.
\end{equation}
As we may assume that $k=\q = 3^{s^w} >5$, we have in particular $w
\geq 1$, and hence $s < 3^{s^w}=\q$. Thus, using
equation~(\ref{E4}), we obtain
\[(\q^{s^{u-w}}-2)(\q^{s^{u-w}}-3)=(q-2)(q-3)  < \q^4  s < \q^5.\]
But, as clearly $u-w \geq 1$ (otherwise, $k=q$, a contradiction to
Corollary~\ref{Cameron_t=5}), this yields a contradiction for $s
\geq 3$. If $s=2$, then $(\q^{2^{u-w}}-2)(\q^{2^{u-w}}-3) < 2\q ^4$
must hold, which cannot be true for $u-w > 1$. Thus, let $u-w=1$.
Hence, it follows from equation~(\ref{E4}) that in particular
\[\q -2 \bigm| (q-2)(q-3)=\q^4-5\q^2+6.\]
But, it is easily seen that $(\q^4-5\q^2+6,\q-2)=(2,\q-2)=1$,
yielding a contradiction.

ad (iv): Let $k=\q$. By condition~(B$^*$), we have
\[2(q-2)(q-3) \left| PSL(2,q)_{0B}\right| = (\q-1)(\q-2)(\q-3)(\q-4)s\]
\[\mbox{with} \;\, \left| PSL(2,q)_{0B}\right| = c \cdot \left\{\begin{array}{ll}
    s,\;\mbox{or}\\
    1.\\
\end{array} \right.\]
As $c \mid \q-1$, we may argue, mutatis mutandis, as in
subcase~(iii). For $k=c\q$, condition~(B$^*$) yields
\[2(q-2)(q-3) \left| PSL(2,q)_{0B}\right| = (c\q-1)(c\q-2)(c\q-3)(c\q-4)s\]
\[\mbox{with} \;\, \left| PSL(2,q)_{0B}\right| = \left\{\begin{array}{ll}
    s,\;\mbox{or}\\
    1.\\
\end{array} \right.\]
We may consider only the equation
\begin{equation}\label{E5}
2(q-2)(q-3)= (c\q-1)(c\q-2)(c\q-3)(c\q-4)s.
\end{equation}
Then, surely $(q-2)(q-3)=q^2-5q+6$ must be divisible by $c\q-3$.
Polynomial division with remainder gives
\begin{eqnarray*}
q^2-5q+6 &=& \bigg(\sum_{i=1}^{m}3^{i-1}\frac{q^2}{(c\q)^i}+
\sum_{j=1}^{\overline{m}} 3^{j-1}
\frac{\big(\big(\frac{3}{c}\big)^m-5 \big)q}{(c\q)^j}\bigg)
\bigg(c\q-3\bigg)\\ & & +\Big(\frac{3}{c}\Big)^{\overline{m}}
\frac{\big(\big(\frac{3}{c})^m-5 \big)q}{\q^{\overline{m}}} +6
\end{eqnarray*}
for a suitable $\overline{m} \in \N$ (such that
\[\mbox{deg}\bigg(\Big(\frac{3}{c}\Big)^{\overline{m}} \frac{\big(\big(\frac{3}{c}\big)^m-5
\big)q}{\q^{\overline{m}}}+6\bigg) < \mbox{deg} \bigg(c\q-3 \bigg)\]
as is well-known). As $c \mid q-1$, clearly $c$ is not divisible by
$3$. Thus, the remainder can be rewritten as
\[\frac{\big(\big(\frac{3}{c}\big)^m-5 \big)}{c^{\overline{m}}} \cdot
3^{s^u- \overline{m}(s^w-1)}+6,\] and hence in order for the
remainder to vanish, necessarily $s^u-\overline{m}(s^w-1)=1$ must
hold. But then, we obtain $3^m=(-2c^{\overline{m}}+5)c^m$, a
contradiction.

ad (v): Let $k=\q+1$. In view of condition~(B$^*$), we have
\[2(q-2)(q-3) \left| PSL(2,q)_{0B}\right| = \q(\q-1)(\q-2)(\q-3)s\]
\[\mbox{with} \;\, \left| PSL(2,q)_{0B}\right| = \frac{\q(\q-1)}{2} \cdot \left\{\begin{array}{ll}
    s,\;\mbox{or}\\
    1.\\
\end{array} \right.\]
Again, it suffices to consider the equation
\begin{equation}\label{E6}
(q-2)(q-3)= (\q-2)(\q-3)s.
\end{equation}
As we may assume that $k=\q +1 = 3^{s^w} +1>5$, it follows in
particular that $w \geq 1$, and hence $s < 3^{s^w}=\q$. Thus, using
equation~(\ref{E6}), we obtain
\[(\q^{s^{u-w}}-2)(\q^{s^{u-w}}-3)=(q-2)(q-3) = (\q - 2)(\q-3)s < ({\q}^2 -2s)(\q-3).\]
But, as clearly $u-w \geq 1$ (otherwise, $k=q+1$, a contradiction to
Corollary~\ref{Cameron_t=5}), this yields a contradiction for every
prime $s$. If $m>1$ even and $k=\q(\q-1)$, then, in view of
condition~(B$^*$), we have
\[2(q-2)(q-3) \left| PSL(2,q)_{0B}\right| = (\q^2-\q-1)(\q^2-\q-2)(\q^2-\q-3)(\q^2-\q-4)s\]
\[\mbox{with} \;\, \left| PSL(2,q)_{0B}\right| = \frac{(\q+1)}{2} \cdot \left\{\begin{array}{ll}
    s,\;\mbox{or}\\
    1.\\
\end{array} \right.\]
We may consider only the equation
\begin{equation*}
(q-2)(q-3)(\q+1)= (\q^2-\q-1)(\q^2-\q-2)(\q^2-\q-3)(\q^2-\q-4)s.
\end{equation*}
As obviously $(\q^2-\q-1,\q+1)=1$, it follows that $\q^2-\q-1 \mid
(q-2)(q-3)$ must hold. But, for $m>1$ even, polynomial division with
remainder gives
\begin{eqnarray*}
q^2-5q+6 &=& \bigg(\sum_{i=1}^{m-1}n_i \frac{q^2}{\q^{i+1}} +
\sum_{j=1}^{m}(n_{j} \cdot n_{m} +
n_{j-1}(n_{m-1}-5))\frac{q}{\q^{j}} \bigg) \bigg(\q^2-\q-1\bigg)\\
& & + (n_{m+1}\cdot
n_{m}+n_m(n_{m-1}-5))\q+n_{m}^2+n_{m-1}(n_{m-1}-5)+6,
\end{eqnarray*}
where $n_i$ denote the $i$-th Fibonacci number recursively defined
via
\[n_0=0,\,\;n_1=n_2=1,\,\; n_{i}=n_{i-1}+ n_{i-2}\,\;(i \geq 3).\]
As it can easily be seen the remainder never vanishes, and hence we
obtain a contradiction. For $k=(\q+1)\q(\q-1)/2$, condition~(B$^*$)
yields
\[2(q-2)(q-3) \left| PSL(2,q)_{0B}\right| = (\frac{\q^3-\q}{2}-1)(\frac{\q^3-\q}{2}-2)(\frac{\q^3-\q}{2}-3)(\frac{\q^3-\q}{2}-4)s\]
\[\mbox{with} \;\, \left| PSL(2,q)_{0B}\right| = \left\{\begin{array}{ll}
    s,\;\mbox{or}\\
    1.\\
\end{array} \right.\]
It suffices to consider the equation
\begin{equation}\label{E8}
2(q-2)(q-3) =
(\frac{\q^3-\q}{2}-1)(\frac{\q^3-\q}{2}-2)(\frac{\q^3-\q}{2}-3)(\frac{\q^3-\q}{2}-4)s.
\end{equation}
If we assume that $\q=3$, then $k=12$. Thus, we obtain from
equation~(\ref{E8}) that $s^u <5$. Hence, there are only a very
small number of possibilities to check, which can easily be ruled
out by hand. Therefore, let us assume that $\q>3$. Then, we have in
particular $w \geq 1$, and hence $s < 3^{s^w}=\q$. Thus, using
equation~(\ref{E8}), we obtain
\[2(q-2)(q-3)< (\frac{\q^3-\q}{2})^4s<\frac{1}{16} \q^{12}s<\frac{1}{16}\q^{13}.\]
On the other hand, it follows that
\begin{eqnarray*}
2(q-2)(q-3)  & = &   (\frac{\q^3-\q}{2}-1)(\frac{\q^3-\q}{2}-2)(\frac{\q^3-\q}{2}-3)(\frac{\q^3-\q}{2}-4)s \\
& \geq & 2(\frac{\q^3-\q}{2}-1)(\frac{\q^3-\q}{2}-2)(\frac{\q^3-\q}{2}-3)(\frac{\q^3-\q}{2}-4)\\
& = & \frac{1}{8}\q^{12}-l
\end{eqnarray*}
with
$l=\frac{1}{2}\q^{10}+\frac{5}{2}\q^9-\frac{3}{4}\q^8-\frac{15}{2}\q^7-17\q^6+\frac{15}{2}\q^5+\frac{279}{8}\q^4
+\frac{95}{2}\q^3-\frac{35}{2}\q^2-50\q-48$. As for $\q>3$, clearly
$l<\frac{1}{16}\q^{12}$ holds, we obtain
\[2(q-2)(q-3) \geq \frac{1}{16}\q^{12}.\]
But as $2(q-2)(q-3)=2(\q^{2m}-5\q^m+6)$, this leaves at most only
$m=6$, which clearly cannot occur since $m=s^{u-w}$.

ad (vi): We may argue, mutatis mutandis, as in subcase~(v).

ad (vii): In view of condition~(B$^*$), we have
\[2(q-2)(q-3) \left| PSL(2,q)_{0B}\right| = (k-1)(k-2)(k-3)(k-4)s\]
\[\mbox{with} \;\, \left| PSL(2,q)_{0B}\right| = \frac{12}{k} \cdot \left\{\begin{array}{ll}
    s,\;\mbox{or}\\
    1.\\
\end{array} \right.\]
It is sufficient to consider the equation
\[(3^{s^u}-2)(3^{s^u}-3)= \frac{k(k-1)(k-2)(k-3)(k-4)}{24}\cdot s.\]
Thus, for $k=6$ respectively $k=12$, we obtain $s^u \leq 2$
respectively $s^u <5$, and thus we have only a very small number of
possibilities to check, which can easily be ruled out by hand.

ad (viii) and (ix): These subcases can be treated similarly to
subcase~(vii), completing the examination of condition~(B$^*$).

\bigskip
\emph{Case} (3): $N=M_v$, $v=11,12,22,23,24$.
\medskip

If $v=12$ or $24$, then $G=M_v$ is always \mbox{$5$-transitive}, and
thus~\cite[Thm.\,3]{Kant1985} yields the designs described in Main
Theorem~\ref{class5-des}. Obviously, flag-transitivity holds as the
$5$-transitivity of $G$ implies that $G_x$ acts block-transitively
on the derived Steiner \mbox{$4$-design} $\D_x$ for any $x \in X$.
By Corollary~\ref{Cameron_t=5}, we obtain for $v=11$ that $k \leq
6$, and for $v=22$ or $23$ that $k \leq 8$, and the very small
number of cases for $k$ can easily be ruled out by hand using
Lemma~\ref{divprop}.

\bigskip
\emph{Case} (4): $N=M_{11}$, $v=12$.
\medskip

As it is known, this exceptional permutation action occurs inside
the Mathieu group $M_{24}$ in its action on $24$ points. This set
can be partitioned into two sets $X_1$ and $X_2$ of $12$ points each
such that the setwise stabilizer of $X_1$ is the Mathieu group
$M_{12}$. The stabilizer in this latter group of a point $x$ in
$X_1$ is isomorphic to $M_{11}$ and operates (apart from its natural
\mbox{$4$-transitive} action on \mbox{$X_1 \setminus \{x\}$})
\mbox{$3$-transitively} on the $12$ points of $X_2$. The geometry
preserved by the $3$-transitive action of $M_{11}$ is not a Steiner
\mbox{$t$-design}, but a $3$-$(12,6,2)$ design
(e.g.~\cite[Ch.\,IV,\,5.3]{BJL1999}).

\medskip

\noindent This completes the proof of Main Theorem~\ref{class5-des}.

\bigskip


\section{The Non-Existence of Flag-transitive Steiner
6-Designs}\label{flag6designs}

We prove the following result:

\begin{mthm}\label{class6-des}
There are no non-trivial Steiner \mbox{$6$-designs} $\D$ admitting a
flag-transitive group \mbox{$G \leq \aut(\D)$} of automorphisms.
\end{mthm}

\bigskip


\subsection{Groups of Automorphisms of Affine
Type}\label{affine typ} \hfill

\bigskip

In the following, we begin with the proof of Main
Theorem~\ref{class6-des}. Using the notation as before, \emph{let us
assume that $\D=(X,\B,I)$ is a non-trivial Steiner \mbox{$6$-design}
with \mbox{$G \leq \aut(\D)$} acting flag-transitively on $\D$
throughout the proof}. Clearly, in the sequel we may assume that
$k>6$ as trivial Steiner \mbox{$6$-designs} are excluded. We will
examine in this subsection successively those cases where $G$ is of
affine type.

\bigskip
\emph{Case} (1): $G \cong AGL(1,8)$, $A \mathit{\Gamma} L(1,8)$, or
$A \mathit{\Gamma} L(1,32)$.
\medskip

We may assume that $k > 6$. If $v=8$, then
Corollary~\ref{Cameron_t=5} would imply that $k=6$. For $v=32$, we
have $\left| G \right| =5v(v-1)$ and Lemma~\ref{divprop} immediately
yields that \mbox{$G \leq \Aut(\D)$} cannot act flag-transitively on
any non-trivial Steiner \mbox{$6$-design} $\D$.

\bigskip
\emph{Case} (2): $G_0 \cong SL(d,2)$, $d \geq 2$.
\medskip

We may argue, mutatis mutandis, as in the corresponding case in Main
Theorem~\ref{class5-des}.

\bigskip
\emph{Case} (3): $G_0 \cong A_7$, $v=2^4$.
\medskip

As $v=2^4$, we have $k \leq 8$ by Corollary~\ref{Cameron_t=5}. But,
Lemma~\ref{Comb_t=5}~(c) obviously eliminates the cases when $k=7$
or $8$.

\bigskip


\subsection{\mbox{Groups of Automorphisms of Almost Simple
Type}}\label{almost simple type} \hfill

\bigskip

We will examine in this subsection successively those cases where
$G$ is of almost simple type.

\bigskip
\emph{Case} (1): $N=A_v$, $v \geq 5$.
\medskip

We may assume that $v \geq 8$. But then $A_v$, and hence also $G$,
is \mbox{$6$-transitive} and does not act on any non-trivial Steiner
\mbox{$6$-design} $\D$ due to~\cite[Thm.\,3]{Kant1985}.

\bigskip
\emph{Case} (2): $N=PSL(2,q)$, $v=q+1$, $q=p^e >3$.
\medskip

For the existence of flag-transitive Steiner \mbox{$6$-designs},
necessarily
\[ r = \frac{q(q-1)(q-2)(q-3)(q-4)}{(k-1)(k-2)(k-3)(k-4)(k-5)} \Bigm|
\left| G_0 \right| \Bigm| \left| P \mathit{\Gamma} L(2,q)_0
\right| = q(q-1)e \] must hold in view of Lemma~\ref{divprop}.
Thus, we have in particular
\begin{equation}\label{Eq-E}
(q-2)(q-3)(q-4) \bigm| (k-1)(k-2)(k-3)(k-4)(k-5)e,\; \mbox{where}\,\; e \leq \mbox{log}_2q.
\end{equation}
But, on the other hand, Corollary~\ref{Cameron_t=5} yields $k \leq
\bigl\lfloor \sqrt{q+1} + \frac{9}{2}\bigr\rfloor <
q^{\frac{1}{2}}+5$. Hence, in view of property~(\ref{Eq-E}), we have
only a small number of possibilities to check, which can easily
be ruled out by hand using Lemma~\ref{Comb_t=5}~(c). Therefore,
\mbox{$G \leq \Aut(\D)$} cannot act flag-transitively on any
non-trivial Steiner \mbox{$6$-design} $\D$. This has also been
proven in~\cite[Cor.\,4.3]{CamPrae1993}, whereas our estimation is slightly better.

\bigskip
\emph{Case} (3): $N=M_v$, $v=11,12,22,23,24$.
\medskip

Due to Corollary~\ref{Cameron_t=5}, we obtain for $v=11$ or $12$
that $k \leq 7$, and for $v=22, 23$ or $24$ that $k \leq 9$, and the very
small number of cases for $k$ can easily be eliminated by hand using
Lemma~\ref{divprop}.

\bigskip
\emph{Case} (4): $N=M_{11}$, $v=12$.
\medskip

By the same arguments as in the corresponding case in Main
Theorem~\ref{class5-des}, it follows that \mbox{$G \leq \Aut(\D)$}
cannot act on any Steiner \mbox{$t$-design} $\D$.

\medskip

\noindent This completes the proof of Main Theorem~\ref{class6-des}.

\bigskip


\subsection*{Acknowledgment}
\mbox{I am grateful to C. Hering for helpful conversations.}
\bibliographystyle{amsplain}
\bibliography{XbibCombDes}
\end{document}